\documentclass[a4paper, twoside,12pt]{article}
\usepackage{fancyhdr}
\usepackage{fnpos}
\usepackage[english]{babel} 
\usepackage[T1]{fontenc}          
\usepackage{graphicx}             
\usepackage{makeidx}
\usepackage{fancybox}
\usepackage{framed}
\usepackage{fancyhdr}
\usepackage[margin=1in]{geometry}
\usepackage{graphicx}
\usepackage{titlesec}
\usepackage{amsmath}
\usepackage{amsfonts}   
\usepackage{amssymb}    
\usepackage{amsthm}
\usepackage{mathtools}
\usepackage{verbatim}   
\usepackage{color}
\usepackage{listings}
\usepackage{graphics}
\usepackage{graphicx}
\usepackage{xcolor}

\newtheorem{theorem}{Theorem}
\newtheorem{proposition}[theorem]{Proposition}
\newtheorem{lemma}[theorem]{Lemma}%
\numberwithin{equation}{section}

\renewcommand{\thefootnote}

\author {Mohamed Gaidi  and Mounir Bedhiafi}
\date{}
\begin{document}
\title{On the Study of the Klein-Gordon Equation in the Dunkl Setting}
\maketitle
\begin{center}
Université Tunis El Manar, Faculté des sciences de Tunis,\\
Laboratoire d’Analyse Mathématique et Applications,\\
LR11ES11, 2092 El Manar I, Tunisie.\\
\textbf{e-mail:} mouhamed.gaidi@esprit.tn, bedhiafi.mounir@yahoo.fr
\end{center}

\begin{abstract}
In Dunkl theory on $\mathbb{R}^{n}$ which generalizes classical Fourier analysis, we study the    solution of the Klein-Gordon-equation defined by:
\begin{eqnarray}
\nonumber \partial_{t}^{2}u-\Delta_{k}u=-m^{2}u \ , \ \ \ u (x,0)=g(x) \ , \ \ \ \partial_{t}u(x,0)=f(x)
\end{eqnarray}
with \ $m > 0$ \ and \ $\partial_{t}^{2}u$ \ is the second derivative of the solution $u$ with respect to $t$ and $\Delta_{k}u$ is the Dunkl Laplacian with respect to $x$ where $f$ and $g$ the two functions in $\mathcal{S}(\mathbb{R}^{n})$ which surround the initial conditions.
We obtain an integral representation for its solution which we gives some properties. As a specific result, we studied the associated energies to the Dunkl-Klein-Gordon equation.
\end{abstract}

\textbf{keyword}\\
Dunkl theory, Klein-Gordon equation, Strichartz estimates, Kinetic and Potential energies

\section{Introduction}\label{sec1}
\par Dunkl operators $T_{i}, \ 1 \leqslant i \leqslant n$  introduced by C.F.Dunkl in \cite{D1}, are differential-difference operators associated with a finite reflection group $G$, acting on some Euclidean space. These operators attached with a positive root system $\mathcal{R}_{+}$ and a non negative multiplicity function $k$, can be considered as perturbations of the usual partial derivatives by reflection parts. During the last years, these operators have gained considerable interest in various fields of mathematics $($see \cite{J1,D1,D2,D3}$)$ and also in physical applications, they are, for example, naturally connected with certain Schr$\ddot{o}$dinger operators and wave equation  \ \cite{AG2,AG,AD,Majj1,Majj2,Per,said1,Sti}. \\
  In a similar of the contexts, we propose applied the Dunkl analysis to the Klein-Gordon equation having many impotance in physics. \\
It was Schr$\ddot{o}$dinger who described the relativistic equation known today as Klein-Gordon, in an attempt to describe the electron within the hydrogen atom.
For a full reconciliation of quantum mechanics with restricted relativity, quantum field theory is needed, in which the Klein-Gordon equation resurfaces as the equation obeyed by the components of all free quantum fields.
In quantum field theory, solutions of the free version which act without interaction always play a role.
However, the Klein-Gordon equation can be written in some versions.
In this context, we will study the Klein-Gordon equation given by:
\begin{eqnarray}\label{Equation}
\partial_{t}^{2}u-\Delta_{x}u = -m^{2}u, \ \ \ u (x,0)=g(x), \ \ \ \partial_{t}u(x,0)=f(x)
\end{eqnarray}
with \ $m > 0$ \ and \ $\partial_{t}^{2}u$ \ is the second derivative of the solution $u$ with respect to $t$ and $\Delta_{x}u$ is the Laplacian with respect to $x$ where $g$ and $f$ the two functions in $\mathcal{S}(\mathbb{R}^{n})$ which surround the initial conditions.
In Peral's article \cite{Per}, he introduce the solution of equation (\ref{Equation}) in the classical case and we have a similar result by applying the Dunkl transform $\mathcal{F}_{k}$ in $x$ to the equation (\ref{Equation}) as:
\begin{eqnarray}\label{Solution}
\mathcal{F}_{k}(u)(\xi,t) = \frac{\sin \big(t (\vert \xi \vert^{2} + m^{2})^{1/2}\big)}{(\vert \xi \vert^{2} + m^{2})^{1/2}} \ \mathcal{F}_{k}(f)(\xi) + \cos \big(t (\vert \xi \vert^{2} + m^{2})^{1/2}\big) \ \mathcal{F}_{k}(g)(\xi)
\end{eqnarray}
Knowing the unicity of the solution of the Dunkl-Klein-Gordon-equation we determine the integral representation for this in respect the mean spherical function of Stein defined by
\begin{eqnarray}\label{Mf}
\mathcal{M}f(x,r) = \frac{1}{d_{k}} \int_{S^{n-1}} \tau_{x}(k) \ f(ry) \ v_{k}(y) \ dw(y), \ \ \ \ \ \big(x \in \mathbb{R}^{n}, \ r > 0\big) 
\end{eqnarray}
where \ \ $d_{k}=\int_{S^{n-1}}v_{k}(x) \ dw(x)$ .
It should be noted that if $m = 0$ in the expression of the operator (\ref{Solution}), we fall directly on the wave equation having, in the Dunkl setting \cite{said1}, for solution:
\begin{eqnarray}
\mathcal{F}_{k}(u)(\xi,t) = \frac{\sin \big(t\vert \xi \vert \big)}{\vert \xi \vert} \ \mathcal{F}_{k}(f)(\xi) + \cos \big(t \vert \xi \vert \big) \ \mathcal{F}_{k}(g)(\xi)
\end{eqnarray}
\par The famous question that arises is the integral representation of the solution of the Klein-Gordon equation (\ref{Solution}) and our interest is to extend the results of Peral \cite{Per} and Strichartz \cite{Sti}. \ On the other hand, Selem Ben Said and B.\O rsted \cite{said1} in their studies of the integral representation of the wave equation considered the locally integrable function on $\mathbb{R}$ defined for $\lambda \in \mathbb{C}$ by 
\begin{eqnarray}
\nonumber \mathbb{S}_{+}^{\lambda}(x) \ = \ x^{\lambda} J_{\lambda}(x) \ \ \ if \ x > 0 \ \ \ and \ \ 0 \ \ if \ \ x \leq 0
\end{eqnarray}
where for $\psi \in D (\mathbb{R})$, the corresponding regular distribution
\begin{eqnarray}
\nonumber \langle \mathbb{S}_{+}^{\lambda}, \psi \rangle \ = \ \int_{0}^{+\infty} x^{\lambda} J_{\lambda}(x) \psi (x) \ dx
\end{eqnarray}
is a holomorphic $D'(\mathbb{R})$- valued function with respect to the variable  $\lambda \in \mathcal{A} \ = \ \{ \lambda \in \mathbb{C} \ / \ \ \Re (\lambda) \ \geq \ 0\}$  where $D(\mathbb{R})$  and $D'(\mathbb{R})$ the space of compactly supported functions and the distributions respectively. 
This last function $\mathbb{S}_{+}^{\lambda}$ plays an important role in the study of the integral representation of the solution (\ref{Solution}) which will be deduced as a function  (\ref{Mf}).
\par The contents of this paper are as follows:\\
In Section $2$, we study the Dunkl theory concerned the Dunkl operator, the Dunkl transform and the Dunkl convolution properties.  We collects some basic definitions and results about harmonic analysis associated with Dunkl operator precisely in the space of tempered distributions. \\
In Section $3$, We study the solution of  Dunkl-Klein-Gordon equation. Indeed, the solution (\ref{Solution}) is deduced directly if we apply the Dunkl transform to the equation (\ref{Equation}).\ 
We show that this last solution (\ref{Solution}) can be written in the form of a convolution product as follows
\begin{eqnarray}
\nonumber \mu_{k}(x,t) \ = \ \big(P_{k,t}^{11} *_{k} f\big)(x) \ + \ \big(P_{k,t}^{12} *_{k} g\big)(x)
\end{eqnarray}
where for a fixed $t$,\ $P_{k,t}^{11}$ and $P_{k,t}^{12}$ are the tempered distributions on $\mathbb{R}^{n}$. 
\par In Section $4$, we deduce two integral representations, the first is with respect to the mean spherical funtion of Stein defined above by (\ref{Mf}) and the second is a similar result with respect to the function:
\begin{eqnarray}\label{Slamda}
 \mathbb{S}_{\lambda}^{+}(x) = \left\{
    \begin{array}{ll}
     x^{\lambda} \ J_{\lambda}(x), & \hbox{ if $x > 0$},\\
0, & \hbox{ if $x \leqslant 0$.}
    \end{array}
  \right.
\end{eqnarray}
For $\psi \in D (\mathbb{R})$, the corresponding regular distribution, is a holomorphic $D' (\mathbb{R})$-valued function with respect to the variable $\lambda \in \mathcal{A} \ = \ \{ \lambda \in \mathbb{C}  \ / \ \ \Re (\lambda) \geqslant  0 \}$. \\
 It admits an analytic continuation into the domain $\mathcal{A}' = \big\{\lambda \in \mathbb{C} \ / \  \lambda \neq 0 , 2 , 3 , ...\big\}$  where with respect to the Dirac distribution $\delta$:
\begin{eqnarray}
\nonumber Res_{\lambda \longmapsto m} x_{+}^{\lambda -1} = \frac{(-1)^{m}}{m!}\delta^{(m)}(x), \ \ \ \ \ for \ \ \ \ \ m=0 , 1 , 2 , 3 , ...
\end{eqnarray}
In second term, we prove a Strichartz estimate satisfied by the solution (\ref{Solution}) deduced  when we extend the variable $t$ on infty.
\par In Section $5$, we prove a similar result obtained in Section $4$ of the total, kinetic and potential energies associated to the Dunkl-Klein-Gordon equation (\ref{Equation}).
\par Along this paper, we denote by:\\
$\bullet \ \ \ \mathcal{E}(\mathbb{R}^{n})$ \ the space of infinitely differentiable functions on $\mathbb{R}^{n}$. \\ 
$\bullet \ \ \ \mathcal{S}(\mathbb{R}^{n})$ the Schwartz space of functions in $\mathcal{E}(\mathbb{R}^{n})$ which are rapidly decreasing as well as their derivatives. \\
$\bullet \ \ \ \mathcal{D}(\mathbb{R}^{n})$ the subspace of $\mathcal{E}(\mathbb{R}^{n})$ of compactly supported functions . \\
$\bullet \ \ \ \mathcal{D}'(\mathbb{R}^{n})$ \ the subspace of distributions defined on $\mathbb{R}^{n}$. \\

\section{Results}\label{sec2}
In this paper, we are interested in the study of the Klein-Gordon equation in the Dunkl setting. \\
We begin by the detemination of its solution by applying the Dunkl transform to the equation (\ref{Equation}) and our result is the following theorem.
\begin{theorem}\label{thhm1}
 The solution $u_{k}$ of the equation (\ref{Equation}) is given by
\begin{eqnarray}\label{Sol}
 u_{k}(\xi,t) = \mathcal{F}_{k}^{-1} \bigg[\frac{\sin \big(t (\vert \xi \vert^{2} + m^{2})^{1/2}\big)}{(\vert \xi \vert^{2} + m^{2})^{1/2}} \ \mathcal{F}_{k}(f)(\xi) + \cos \big(t (\vert \xi \vert^{2} + m^{2})^{1/2}\big) \ \mathcal{F}_{k}(g)(\xi)\bigg]
\end{eqnarray}
\end{theorem}
Consequently, we can write the operator (\ref{Sol}) as a convolution product as the following theorem.
\begin{theorem}\label{thhm2}
The solution (\ref{Sol}) of the Cauchy problem (\ref{Equation}) is given by
\begin{eqnarray}
u_{k}(x,t) = \big(P_{k,t}^{11} *_{k} f\big)(x) + \big(P_{k,t}^{12} *_{k} g\big) (x),
\end{eqnarray}
where for a fixed  $t, \ P_{k,t}^{11}$  and  $P_{k,t}^{12}$  are the tempered distributions on $\mathbb{R}^{n}$ given by 
\begin{eqnarray}\label{distrb}
 \nonumber P_{k,t}^{11} &=& \mathcal{F}_{k}^{-1} \big[\cos \big(t (\vert . \vert^{2} + m^{2})^{1/2}\big)\big] \\
 P_{k,t}^{12} &=& \mathcal{F}_{k}^{-1} \big[\frac{\sin \big(t (\vert . \vert^{2} + m^{2})^{1/2}\big)}{(\Vert . \Vert^{2} + m^{2})^{1/2}}\big] 
\end{eqnarray}
\end{theorem}
In second term, we determine the integral representation of the operator (\ref{Sol}) with respect to the mean spherical function Stein $\mathcal{M}f$ defined by (\ref{Mf}) and our main result is the following theorem.
\begin{theorem}\label{thp}
Let  $u_{k}$ the solution the Dunkl Klein-Gordon equation (\ref{Equation}). For all  $(x,t) \in \mathbb{R}^{n} \times \mathbb{R}$, we have
\begin{eqnarray}
\nonumber u_{k}(x,t) &=& \mathcal{L}_{n , k , m} \  \int_{0}^{\vert t \vert} r^{2 \gamma_{k}+n-1} \   \mathbb{S}_{-\gamma_{k}-n/2 +1/2} [m(t^{2}-r^{2})] \ \mathcal{M}_{f}(x,r) \ dr \\
\nonumber &+& \mathcal{L}_{n , k , m} \ \int_{0}^{\vert t \vert} r^{2 \gamma_{k}+n-1} \ \frac{d}{dt} \big( \mathbb{S}_{-\gamma_{k}-n/2+1/2} [m(t^{2}-r^{2})] \big) \ \mathcal{M}_{g}(x,r) \ dr 
\end{eqnarray}
where $\mathcal{M}f$ and \ $\mathbb{S}_{\lambda,\ (\lambda \in  \{\lambda \in \mathbb{C} / \ \ \Re (\lambda) \geqslant 0 \})}$ are these functions defined above by (\ref{Mf}) and (\ref{Slamda}) successively and 
$$\mathcal{L}_{n , k , m} \ = \ d_{k} \frac{\sqrt{\pi}}{\Gamma (\gamma_{k}+n/2)} \ 2^{-\gamma_{k}-n/2+1/2} \ m^{2 \gamma_{k}+n-1}$$ 
\end{theorem}
By writing the $L^{2}$-norm of the solution $u_{k}$ we can prove that as $\vert t \vert \longrightarrow + \infty$, the function $t \longmapsto \Vert u_{k}(.,t) \Vert_{2}$ has a finite limit depending on the initial datas. Consequently, we obtain the following lemma.
\begin{lemma}\label{ESTM}
 Let $u_{k}$ the solution of the Klein-Gordon equation given by (\ref{Equation}).\\
 As \ $\vert t \vert \longmapsto + \infty$, \ the function   $t \longmapsto \Vert u_{k} (. , t) \Vert_{k}$  has a finite limit depending on the initial datas as:
 \begin{eqnarray}\label{estm}
 \lim_{\vert t \vert \longmapsto + \infty } \Vert u_{k} (. , t) \Vert_{k}^{2} \ \leqslant \ \frac{1}{2} \ \Vert (- \Delta_{k})^{-1/2} f \Vert_{k}^{2} + \frac{1}{2} \ \Vert g\Vert_{k}^{2} 
 \end{eqnarray}
 Here  $\Vert.\Vert_{k}$ \ denotes the norm in  $L^{2} (\mathbb{R}^{n} ; w_{k}(x) dx)$.
\end{lemma}
It's remembered that the last Strichartz  type inequality (\ref{ESTM}) have a big importance in physical. \\
After the above discussions, we can study the energies associated to Klein-Gordon equation (\ref{Equation}). \\
Indeed, the total enargy of $u_{k}$ is given by
\begin{eqnarray}
\nonumber \mathcal{E}_{k}[u_{k}](t) = \mathcal{K}_{k}[u_{k}](t) + \mathcal{P}_{k}[u_{k}](t) 
\end{eqnarray}
where $\mathcal{E}_{k}, \mathcal{K}_{k}$ and $\mathcal{P}_{k}$ are the total, kinetic and potential energies successively. \\
Using some trigonometric formulas, we shall determine the expression of $\mathcal{K}_{k}, \mathcal{P}_{k}$ and $\mathcal{E}_{k}$ successively and using the Riemann-Lebesgue theorem, we calculate the limits of these three energies as $\vert t \vert \longrightarrow + \infty$ and our goal is to prove the following proposition.
\\
\begin{proposition}\label{collect}
Let $\mathcal{K}_{k}, \mathcal{P}_{k}$ and $\mathcal{E}_{k}$ the kinetic, potential and total energies successively associated to the solution $u_{k}$ of the Dunkl Klein-Gordon equation (\ref{Equation}). \\
Then:
\begin{eqnarray}
\nonumber \lim_{\vert t \vert \longmapsto + \infty} \mathcal{K}_{k}[u_{k}](t) &=&  \frac{c_{k}^{-2}}{4} \Vert f \Vert_{k}^{2} \ + \  \frac{c_{k}^{-2}}{4} \Vert T_{j}^{k}(g) \Vert_{k}^{2} \ + \  \frac{c_{k}^{-2}}{4} m^{2} \Vert g \Vert_{k}^{2} \\
\nonumber  \lim_{\vert t \vert \longmapsto + \infty} \mathcal{P}_{k}[u_{k}](t) &=& \frac{c_{k}^{-2}}{4} \Vert f \Vert_{k}^{2} \ + \ \frac{c_{k}^{-2}}{4} \Vert T_{j}^{k}(g) \Vert_{k}^{2} - \frac{c_{k}^{-2}}{4} m^{2} \Vert (-\Delta_{k}+m^{2})^{-1/2} f \Vert_{k}^{2} \\
\nonumber \lim_{\vert t \vert \longmapsto + \infty} \mathcal{E}_{k}[u_{k}](t) \ &=& \ \frac{c_{k}^{-2}}{2} \Vert f \Vert_{k}^{2} \ + \ \frac{c_{k}^{-2}}{2} \Vert T_{j}^{k}(g) \Vert_{k}^{2} \ + \ \frac{m^{2} \ c_{k}^{-2}}{4} \Vert g \Vert_{k}^{2} \\
\nonumber &-& \frac{c_{k}^{-2}}{4} m^{2} \Vert (-\Delta_{k}+m^{2})^{-1/2} f \Vert_{k}^{2}
\end{eqnarray}
\end{proposition}
where $\Delta_{k}$ is the Dunkl-Laplacian operator given by (\ref{lap}).

\section{On the Dunkl Theory}\label{sec3}
\par In this section, we recall some results on Dunkl theory ( see \cite{J1,D1,D2,D3} and we refer for more details to the surveys \cite{Tr1}.
 \par Let \ $G \subset O (\mathbb{R}^{n})$ \ be a finite reflection group on \ $\mathbb{R}^{n}$, associated with a root system $R$. For $\alpha \in R$, we denote by $\mathbb{H}_{\alpha}$ \ the hyperplane orthogonal to $\alpha$. For a given \ $\beta \in \mathbb{R}^{n}\setminus \cup_{\alpha \in R}\mathbb{H}_{\alpha}$. We fix a positive subsystem $R_{+} = \{\alpha \in R: \ \langle \alpha, \beta\rangle > 0\}$. We denote by $k$ a nonnegative multiplicity function defined on $R$ with the property that $k$ is G-invariant. We associate with $k$ the index 
 \begin{eqnarray}
 \nonumber \gamma = \sum_{\xi \in R_{+}} k (\xi),
 \end{eqnarray}
 and a weighted measure $\nu_{k}$ given by
 \begin{eqnarray}
 \nonumber d\nu_{k}(x) = w_{k}(x) dx \ \ \ \ where \ \ \ \ w_{k}(x) = \Pi_{\xi \in R_{+}} \vert\langle\xi , x\rangle\vert^{2 k (\xi)}, \ \ \ \ x \in \mathbb{R}^{n}, 
 \end{eqnarray}
 Further, we introduce the Mehta-type constant $c_{k}$ by
 \begin{eqnarray}
 \nonumber c_{k} = \bigg(\int_{\mathbb{R}^{n}} e^{-\frac{\Vert x \Vert^{2}}{2}} w_{k}(x) dx\bigg)^{-1}
 \end{eqnarray}
 For every $1 \leqslant p \leqslant +\infty$, we denote by $L_{k}^{p}(\mathbb{R}^{n})$, the spaces 
 $L^{p}(\mathbb{R}^{n} , d \nu_k (x))$. \\
 \par For an integrable function $f$ on $\mathbb{R}^{n}$ with respect to the measure $w_{k}(x)dx$ we have the relation
 \begin{eqnarray}\label{hom}
 \int_{\mathbb{R}^{n}}f(x) w_{k}(x) dx \ = \ \int_{0}^{+\infty} \bigg( \int_{\mathcal{S}^{n-1}}f(r \beta) w_{k}(r\beta) d\sigma(\beta) \bigg) r^{n-1} dr ,
 \end{eqnarray}
 where $d\sigma$ is the normalized surface measure on the unit sphere $\mathcal{S}^{n-1}$ on $\mathbb{R}^{n}$. \\
 Using the homogeneity of $w_{k}$, the relation (\ref{hom}) can also be written in the form
 \begin{eqnarray}\label{hoj}
 \int_{\mathbb{R}^{n}}f(x) w_{k}(x) dx \ = \ \int_{0}^{+\infty} \bigg(\int_{\mathcal{S}^{n-1}} f(r \beta) w_{k}(\beta) \bigg) r^{2\gamma_{k}+n-1} dr
 \end{eqnarray}
 In particular if $f$ is radial $($ i . e . $SO(n)$-invariant $)$, then there exists a function $F$ on $[0 , +\infty)$ such that  $f(x) = F(\Vert x \Vert) = F(r)$ with $r = \Vert x \Vert$, and the relation (\ref{hoj}) takes the form 
 \begin{eqnarray}\label{mes}
  \int_{\mathbb{R}^{n}} f(x) \ d\nu_{k}(x) = d_{k} \ \int_{0}^{+\infty} F(r) \ r^{2\gamma_{k}+n-1} dr,
 \end{eqnarray}
 
 \par The Dunkl operators $T_{j}, \ 1 \leqslant j \leqslant n$, \ on $\mathbb{R}^{n}$ associated with the reflection group $G$ and the multiplicity function $k$ are the first-order differential-difference operators given by
 \begin{eqnarray}
 \nonumber T_{j}f(x) = \frac{\partial f}{\partial x_{j}} + \sum_{\alpha \in R_{+}} k(\alpha) \ \alpha_{j} \ \frac{f(x) - f(\rho_{\alpha} (x))}{\langle \alpha, x \rangle}, \ \ \ \ \ f \in \mathcal{E}(\mathbb{R}^{n}), \ \ x \in \mathbb{R}^{n}
 \end{eqnarray}
 where $\rho_{\alpha}$ is the reflection on the hyperplane $\mathbb{H}_{\alpha}$ and $\alpha_{j}=\langle \alpha, e_{j}\rangle, \ (e_{1} , ... , e_{n})$ \ being the canonical basis of $\mathbb{R}^{n}$. \\
 Related to the Dunkl operator, the Dunkl Laplacian operator is defined by
 \begin{eqnarray}\label{lap}
 \Delta_{k} \ = \ \sum_{i=1}^{n} ( T_{j} )^{2}
 \end{eqnarray}
 In the case $k \equiv 0$, the wheighed function $w_{k} \equiv 1$ and the measure $\nu_{k}$ associated to the Dunkl operators coincide with the Leabesgue measure. The $T_{j}$ reduce to the corresponding partial derivatives. Therefore Dunkl analysis can be viewed as a generalization of classical Fourier analysis.
 \par For $y \in \mathbb{C}^{n}$, the system
 \begin{eqnarray}
\nonumber T_{\xi}^{k} f \ = \ \langle y , \xi \rangle \ f , \ \ \ \ for \ all \ \xi \in \mathbb{R}^{n}
 \end{eqnarray}
 admits a unique analytic solution on $\mathbb{R}^{n}$, denoted by $E_{k}(. , y)$ \ and called the Dunkl kernel. This kernel has a unique holomorphic extension to $\mathbb{C}^{n}\times \mathbb{C}^{n}$. 
 \par The Dunkl transform $\mathcal{F}_{k}$ is defined for $f \in \mathcal{D}(\mathbb{R}^{n})$ by
 \begin{eqnarray}\label{dnk}
  \mathcal{F}_{k}(f)(x) = c_{k} \int_{\mathbb{R}^{n}} f(y) \ E_{k}(-ix , y) \ d\nu_{k}(y) , \ \ \ \ x \in \mathbb{R}^{n}
 \end{eqnarray}
  It is remembered that if $f$ a radial function in $L^{1}(\mathbb{R}^{n} , w_{k}(x) dx)$ such that \ $f(x) = \tilde{f}(\vert x \vert ), \ x \in \mathbb{R}^{n}$. Then $\mathcal{F}_{k}(f)$ is also radial and
 \begin{eqnarray}\label{rad}
 \mathcal{F}_{k}(f)(x) \ = \ d_{k} \int_{0}^{+\infty} \tilde{f}(s) \ \mathcal{J}_{\gamma_{k}+n/2-1}(s \vert x \vert ) \ s^{\gamma_{k}+n/2} \ ds , \ \ \ x \in \mathbb{R}^{n}
 \end{eqnarray}
 where \ \ $d_{k} \ = \ \vert x \vert^{-(\gamma_{k}+n/2-1)}/ w_{n}$ and $\mathcal{J}_{\alpha}$ is the normalized Bessel function of indice $\alpha$.\\
 The Dunkl transform satisfies: \\
 $(i)$ \ \ The Dunkl transform of a function $f \in L_{k}^{1}(\mathbb{R}^{n})$ has the following basic property
 \begin{eqnarray}
 \nonumber \Vert \mathcal{F}_{k}(f)\Vert_{\infty , k}  \ \leqslant \ \Vert f \Vert_{1,k}
 \end{eqnarray}
 $(ii)$ \ \ The Dunkl transform is an automorphism on the Schwartz space $\mathcal{S}(\mathbb{R}^{n})$. \\
 $(iii)$ \ \ When both $f$ and $\mathcal{F}_{k}(f)$ are in $L_{k}^{1}(\mathbb{R}^{n})$, we have the inversion formula
 \begin{eqnarray}
 \nonumber f(x) = \int_{\mathbb{R}^{n}} \mathcal{F}_{k}(f)(y) \ E_{k}(ix , y) \ d\nu_{k}(y) , \ \ \ \ x \in \mathbb{R}^{n}
 \end{eqnarray}
 The Dunkl translation operator $\tau_{x}, \ x \in \mathbb{R}^{n}$, was introduced in \cite{Tr2} on $\mathcal{E}(\mathbb{R}^{n})$. \\
 For $f \in \mathcal{S}(\mathbb{R}^{n})$, \ we have
 \begin{eqnarray}\label{kern}
  \mathcal{F}_{k}(\tau_{x}(f))(y) = E_{k}(ix , y) \ \mathcal{F}_{k}(f)(y),
 \end{eqnarray}
In the case when $f(x) = \tilde{f}(\vert x \vert )$ is a continuous radial function that belongs to $L^{2}(\mathbb{R}^{n}, w_{k}(x) dx)$, the Dunkl translation is represented by
\begin{eqnarray}
\nonumber \tau_{x}(f)(y) = \int_{\mathbb{R}^{n}}\tilde{f} (\sqrt{\vert y \vert^{2} + \vert x \vert^{2} + 2 \langle y, \eta \rangle}) \ d\nu_{x}(\eta)
\end{eqnarray}
This formula shows that the Dunkl translation operators can be extended to all radial functions $f$ in $L^{p}(\mathbb{R}^{n}, w_{k}(x) dx), \ 1 \leq p \leq \infty $ and the following holds
\begin{eqnarray}
\nonumber \Vert \tau_{x} (f)\Vert_{p,k} \ \leqslant \Vert f \Vert_{p,k}
\end{eqnarray}
We define the convolution product for suitable funcions $f$ and $g$ by
\begin{eqnarray}\label{yi}
 f *_{k} g (x) = \int_{\mathbb{R}^{n}} \tau_{x} (f)(-y) \ g(y) \ d\mu_{k}(y), \ \ \ \ \ x \in \mathbb{R}^{n}
\end{eqnarray}
We note that it is commutative and satisfies the following property
\begin{eqnarray}\label{for}
 \mathcal{F}_{k} (f *_{k} g ) \ = \ \mathcal{F}_{k}(f) \ \mathcal{F}_{k}(g)
\end{eqnarray}
\par Next we turn our attention to the Dunkl convolution of two distributions. Denote by $D(\mathbb{R}^{n})$ the space of smoothly compact supported functions on $\mathbb{R}^{n}$ and set $D'(\mathbb{R}^{n})$ 
to be its dual.\\ 
Let $f \in L^{1}(\mathbb{R}^{n} , w_{k}(x) dx)$ and $\phi \in D(\mathbb{R}^{n})$. Let $\tau_{f}$ to be the linear form on $D(\mathbb{R}^{n})$ defined by:
\begin{eqnarray}
\nonumber \big<\tau_{f} , \phi\big> = \int_{\mathbb{R}^{n}} f(x) \ \phi(x) \ w_{k}(x) \ dx
\end{eqnarray}
We used the fact that

\begin{eqnarray}
\nonumber \int_{\mathbb{R}^{n}} \phi (y) \ \tau_{x}(k) \ \psi (-y) \ w_{k}(y) \ dy = \int_{\mathbb{R}^{n}} \tau_{y}(k) \ \phi (x) \ \psi (y) \ w_{k}(y) \ dy. 
\end{eqnarray}
Therefore
\begin{eqnarray}
\nonumber S*_{k}\phi (x) = \big< S , \tau_{x}(k) \phi\big>.
\end{eqnarray}
Clearly
\begin{eqnarray}
\nonumber S*_{k}\phi = \phi *_{k} S.
\end{eqnarray}
Since the mapping $\phi \longrightarrow \mathcal{F}_{k}(\phi)$ of $\mathcal{S}(\mathbb{R}^{n})$ onto $\mathcal{S}(\mathbb{R}^{n})$ is linear and continuous in the topology of $\mathcal{S}(\mathbb{R}^{n})$, we can now define the Dunkl transform of a tempered distribution $\tau$ as the tempered distribution $\mathcal{F}_{k}(\tau)$ defined through,
\begin{eqnarray}
\nonumber \big<\mathcal{F}_{k}(\tau) , \phi \big> = \big< \tau , \mathcal{F}_{k} (\phi) \big>, \ \ \ \ \ \ \ \phi \in \mathcal{S}(\mathbb{R}^{n})
\end{eqnarray}
Notice that, if $\tau \in \mathcal{E}'(\mathbb{R}^{n})$, \ then we  can write
\begin{eqnarray}
\nonumber \big<\mathcal{F}_{k}(\tau) , \phi \big> = \big< \big< \tau_{y} , E_{k}(-ix , y) \big> \phi \big>.
\end{eqnarray}
That is
\begin{eqnarray}
\nonumber \mathcal{F}_{k} (\tau) (\xi) = \big<\tau_{x} , E_{k} (-i \xi , x )\big>, \ \ \ \ \ \ \ \forall \tau \in \mathcal{E}'(\mathbb{R}^{n}).
\end{eqnarray}
\par Next we turn our attention to the behavior of the convolution of two distribution under the Dunkl transform. We claim that:
\begin{eqnarray}
\nonumber \mathcal{F}_{k} \big(S*_{k}\tau\big) = \mathcal{F}_{k} (S) \ . \ \mathcal{F}_{k} (\tau)
\end{eqnarray}
holds if one of the distribution is of compact support and the other one is a tempered distribution. \\
On the other hand, since \ $S \in \mathcal{E}'(\mathbb{R}^{n}), \ \mathcal{F}_{k}(S)$ \ is a $\mathcal{C}^{\infty}$  slowly increasing function and
\begin{eqnarray}
\nonumber \chi (y) &=& \int_{\mathbb{R}^{n}} \mathcal{F}_{k}(S)(\xi) \ E_{k}(-iy , i \xi) \ \phi (\xi) \ w_{k}(\xi) \ d \xi \\
\nonumber &=& \mathcal{F}_{k} \big(\mathcal{F}_{k}(S). \phi\big)(y)
\end{eqnarray}
Thus, $\chi (y)$ is in $\mathbb{S}(\mathbb{R}^{n})$, and the mapping \ $\mathcal{S}(\mathbb{R}^{n}) \longrightarrow \mathcal{S}(\mathbb{R}^{n})$, defined by \ $\phi \longmapsto \chi$ is continuous, also, if $\tau$ is a tempered distribution, then
\begin{eqnarray}
\nonumber \big<\mathcal{F}_{k}(S*_{k}\tau), \phi \big> = \big<\tau, \chi \big>
\end{eqnarray}
is well defined, depends continuously on $\phi$ and 
\begin{eqnarray}
\nonumber \big<\mathcal{F}_{k}(S*_{k}\tau), \phi \big> &=& \big< \tau, \mathcal{F}_{k} (\mathcal{F}_{k}(S). \phi)\big> \\
\nonumber &=& \big< \mathcal{F}_{k} (\tau), \mathcal{F}_{k}(S). \phi \big> \\
\nonumber  &=& \big<\mathcal{F}_{k} (S) \ \mathcal{F}_{k}(\tau), \phi\big> 
\end{eqnarray}
That is
\begin{eqnarray}\label{formula}
 \mathcal{F}_{k}(S*_{k}\tau) = \mathcal{F}_{k} (S) \ \mathcal{F}_{k}(\tau) , \ \ \ \ \ for \ all \ \ \ \ S \in \mathcal{E}'(\mathbb{R}^{n}), \ \tau \in \mathcal{S}'(\mathbb{R}^{n})
\end{eqnarray}

\section{Solution \ of \ the \ Dunkl \ Klein \ - \ Gordon \ equation}
Except in a few places, most of the results below hold for complex-valued multiplicity functions $k$ such that $\Re (k) \geqslant 0$.\\
Henceforth, $\mathcal{K}^{+}$ denotes the set of multiplicity functions $k = (k_{\alpha})_{\alpha \in \mathcal{R}}$ such that $k_{\alpha}\geqslant 0$ for all $\alpha \in \mathcal{R}$.
However, for the reader's convenience, we will restrict ourselves to multipkicity functions $k \in \mathcal{K}^{+}$.
\par For $k$ in $\mathcal{K}^{+}$, consider  the Klein-Gordon equation associated with the Dunkl-Laplacian operator defined by 
\begin{eqnarray}\label{Klein}
 \partial_{t}^{2}u-\Delta_{x}u=-m^{2}u \ , \ \ \ u (x,0)=g(x) \ , \ \ \ \partial_{t}u(x,0)=f(x)
\end{eqnarray}
with $m > 0$ and \ $\partial_{t}^{2}u$ \  is the second derivative of the solution $u$ with respect to $t$ and $\Delta_{x}u$ is the Laplacian with respect to $x$ and $g$ and $f$ the two functions in $\mathcal{S}(\mathbb{R}^{n})$ which surround the initial conditions. \\
Here the functions $f$ and $g$ belongs to $\mathcal{S}(\mathbb{R}^{n})$. The subscript  $t$ indicates differentiation in the $t-$variable. 
\par Next, we will prove that the solution $u_{k}$ is expressed as a sum of \ $*_{k}$ - \  convolutions of $f$ and $g$, taking Dunkl transform in the $x$ variables and using the formula (\ref{for}) we get the solution of the Klein-Gordon equation given by the operator (\ref{Sol}). 
\par Now,  we only assume $k \in \mathcal{K}^{+}$  and $N \geqslant 1$. For $t \in \mathbb{R}$, denote by $P_{k,t}$ the $2\times 2$ matrix of tempered distributions on $\mathbb{R}^{n}$,
$$
P_{k,t} = 
\begin{bmatrix}
P_{k,t}^{11}&P_{k,t}^{12}\\
P_{k,t}^{21}&P_{k,t}^{22}\\
\end{bmatrix}
$$
where
\begin{eqnarray}
\nonumber P_{k,t}^{11} &=& P_{k,t}^{22} = \mathcal{F}_{k}^{-1} \big[\cos \big(t (\vert \xi \vert^{2} + m^{2})^{1/2}\big)\big], \\
\nonumber P_{k,t}^{12} &=& \mathcal{F}_{k}^{-1} \big[\frac{\sin \big(t (\vert \xi \vert^{2} + m^{2})^{1/2}\big)}{(\vert \xi \vert^{2} + m^{2})^{1/2}}\big], \ \ \ \ \ \ and \\
\nonumber P_{k,t}^{21} &=& \mathcal{F}_{k}^{-1} \big[-(\vert \xi \vert^{2} + m^{2})^{1/2} \ \sin \big(t (\vert \xi \vert^{2} + m^{2})^{1/2}\big)\big]
\end{eqnarray}
Put \ $U_{k}(x,0)=\begin{bmatrix}
f(x)\\
g(x)\\
\end{bmatrix}$, \ where the Cauchy data \ $(f,g) \in \mathcal{S}(\mathbb{R}^{n}) \times \mathcal{S}(\mathbb{R}^{n})$. \\
Thus, we may define the vector column $U_{k}(x,t)$ by 
\begin{eqnarray}\label{colum}
\nonumber U_{k}(x,t) &=& \big\{P_{k,t} *_{k} U_{k}(. , 0)\big\} (x) \\
 &=& \begin{bmatrix}
P_{k,t}^{11}&P_{k,t}^{12}\\
P_{k,t}^{21}&P_{k,t}^{22}\\
\end{bmatrix}
*_{k}\begin{bmatrix}
f(x)\\
g(x)\\
\end{bmatrix}
(x) 
\end{eqnarray}
By applying the Dunkl transform $\mathcal{F}_{k}$ to (\ref{colum}), in the $x-$variable, we get:
\begin{eqnarray}
\nonumber \mathcal{F}_{k}\big(U_{k}(.,t)\big) (\xi) = e^{t \mathbb{A}} \ \mathcal{F}_{k}\big(U_{k}(.,0)\big) (\xi),
\end{eqnarray}
where 
\begin{eqnarray}
\nonumber \mathbb{A} = \begin{bmatrix}
0 & 1\\
- \vert \xi \vert^{2} - m^{2} & 0\\
\end{bmatrix}
\end{eqnarray}
That is \ $\mathcal{F}_{k}\big(U_{k}(.,t)\big)(\xi)$ \ is a solution to the following ordinary differential equation:
\begin{eqnarray}\label{ord}
\nonumber \partial_{t} \mathcal{F}_{k}\big(U_{k}(.,t)\big)(\xi) &=& \mathbb{A} \ \mathcal{F}_{k}\big(U_{k}(.,t)\big)(\xi) \\
&=& \begin{bmatrix}
0 & 1\\
- \vert \xi \vert^{2} - m^{2} & 0\\
\end{bmatrix} 
\mathcal{F}_{k}\big(U_{k}(.,t)\big)(\xi)
\end{eqnarray}
Using the fact that \ $-\vert \xi \vert^{2} \mathcal{F}_{k} (f)(\xi) =  \mathcal{F}_{k} (\Delta_{k}f)(\xi)$  \ and the injectivity of the Dunkl transform, we deduce on the basing of  (\ref{ord}):
\begin{eqnarray}
\nonumber \partial_{t}U_{k}(x,t) = \begin{bmatrix}
0 & 1\\
\Delta_{k} - m^{2} & 0\\
\end{bmatrix} 
U_{k}(x,t)
\end{eqnarray}
Thus, if we write \ $U_{k}(x,t) = \begin{bmatrix}
u_{k}(x,t) \\
v_{k}(x,t) \\
\end{bmatrix},
$
then \ $u_{k}(x,t)$ \ satisfies the following Klein-Gordon equation 
\begin{eqnarray}
\nonumber \partial_{t}^{2} u_{k}(x,t)- \Delta_{k} u_{k}(x,t) = -m^{2} \ u_{k}(x,t) 
\end{eqnarray}
Moreover, from \ (\ref{colum}) \ and in the light of very last fact pointed out in the previous section regarding \ $*_{k}, \ u_{k}(x,t) \in \mathcal{S}(\mathbb{R}^{n})$ for each \ $t \in \mathbb{R}$. 
\par Furthermore  \ $u_{k}(x,t) \longrightarrow f(x)$ \ as \ $t \longmapsto 0$. \\
Indeed,  if \ $\delta$ \ denotes the Dirac functional, then, as \ $t \longmapsto 0$ \ we have \  $\mathcal{F}_{k}^{-1} \big(\cos \big(t (\vert \xi \vert^{2} + m^{2})^{1/2}\big) \big) \longrightarrow  \ \delta$ 
\ in $\mathcal{S}^{'} (\mathbb{R}^{n})$ and thus in $\mathcal{D}^{'}(\mathbb{R}^{n})$. \\
On the other hand:
\begin{eqnarray}
\nonumber \mathcal{F}_{k}^{-1} \bigg(\frac{\sin \big(t (\vert \xi \vert^{2} + m^{2})^{1/2}\big)}{(\vert \xi \vert^{2} + m^{2})^{1/2}}\bigg) \longrightarrow 0 \ \ \ as \ \ \ t \longmapsto 0
\end{eqnarray}
Using the continuity of the Dunkl convolution \ $*_{k}$, we deduce that 
\begin{eqnarray}
\nonumber u_{k}(x,t) \longrightarrow \big(\delta *_{k} f\big) (x) = f(x), \ \ \ \ as \ \ \ \ t \longmapsto 0
\end{eqnarray}
Similary, one can prove that 
\begin{eqnarray}\label{partial}
 \big(\partial_{t} u_{k} \big) (x,t) \longrightarrow g(x) \ \ \ \ \ as \ \ \ \ \ t \longmapsto 0
\end{eqnarray}
The estimate (\ref{partial}) \ is deduced directly  if we write
\begin{eqnarray}
\nonumber \partial_{t}u_{k}(x,t) &=& \mathcal{F}_{k}^{-1}\bigg(\cos \big(t (\Vert \xi \Vert^{2}+m^{2})^{1/2}\big) *_{k} g\bigg)(x) \\
\nonumber &+& \mathcal{F}_{k}^{-1}\bigg( (\Vert \xi \Vert^{2}+m^{2})^{1/2} \ \sin \big(t (\Vert \xi \Vert^{2}+m^{2})^{1/2}\big) *_{k} f\bigg)(x).
\end{eqnarray}
This completes the proof of theorem \ref{thhm1}. \\
Now, we collects all the above facts and discussions to write the solution (\ref{Sol}) of the Dunkl-Klein-Gordon equation as convolution product. Indeed, using the formula (\ref{formula}), the two tempered distributions (\ref{distrb}) holds directly and this gives the theorem \ref{thhm2}.

\section{Integral \ representation \ of \ the \ solution \ of \ the \ Dunkl  \ Klein \ - \ Gordon \ equation}\label{sec4}

\par Let $\mathcal{A} = \big\{ \lambda \in \mathbb{C} / \ \Re (\lambda ) \geqslant 0 \big\}$. Consider the locally integrable function on $\mathbb{R}$ defined for $\lambda \in \mathbb{C}$ by 
\begin{eqnarray}\label{Slmd}
 \mathbb{S}_{+}^{\lambda} = \left\{
    \begin{array}{ll}
     x^{\lambda} \ J_{\lambda}(x), & \hbox{ if $x > 0$},\\
0, & \hbox{ if $x \leqslant 0$.}
    \end{array}
  \right.
\end{eqnarray}
For $\psi \in D (\mathbb{R})$, the corresponding regular distribution 
\begin{eqnarray}
\nonumber \langle  \mathbb{S}_{+}^{\lambda} , \psi \rangle \ = \ \int_{0}^{+\infty} x^{\lambda} \ J_{\lambda}(x) \ \psi (x) \ dx,
\end{eqnarray}
is a holomorphic $D'(\mathbb{R})-$ valued function with respect to the variable $\lambda \in \mathcal{A}$. 
\par Next, we turn our attention to the relation between $u_{k}$ and the spherical mean operator  $f \longmapsto \mathcal{M}f$  defined in \cite{said1} for  $\mathcal{C}^{\infty}(\mathbb{R}^{n})$  by
\begin{eqnarray}\label{MF}
\mathcal{M}f(x,r) = \frac{1}{d_{k}} \int_{S^{n-1}} \tau_{x}(k) \ f(ry) \ v_{k}(y) \ y, \ \ \ \ \ \big(x \in \mathbb{R}^{n}, \ r \geqslant 0\big) 
\end{eqnarray}
where \  $d_{k}=\int_{S^{n-1}}v_{k}(x) \ dw(x)$  and  $S^{n-1}$ \ is the unit sphere on $\mathbb{R}^{n}$.
 By theorem \ref{thhm2} and the identity (\ref{kern}), we know that
\begin{eqnarray}\label{Uk}
 u_{k}(x,t) &=& \int_{\mathbb{R}^{n}} P_{k,t}^{11}(y) \ \tau_{x}(k) \ f(y) \ v_{k}(y) \ dy \\
 &+& \int_{\mathbb{R}^{n}} P_{k,t}^{12}(y) \ \tau_{x}(k) \ g(y) \ v_{k}(y) \ dy
\end{eqnarray}
 By R$\ddot{o}$sler and Voit \cite{Voit}, if \  $F(x) = F_{0}(\Vert \xi \Vert)$  where  \ $F_{0} : \ \mathbb{R}_{+} \longrightarrow \mathbb{C}$, then
\begin{eqnarray}\label{F01}
\mathcal{F}_{k} \ F(\xi) = \mathcal{H}_{\gamma_{k}+n/2-1} \ F_{0}(\Vert \xi \Vert),
\end{eqnarray}
where \ $\mathcal{H}_{\alpha}$ \ denotes the Hankel transform defined by
\begin{eqnarray}\label{F02}
\mathcal{H}_{\alpha} F_{0}(r) = \frac{1}{2^{\alpha} \Gamma (\alpha +1)} \int_{0}^{+\infty} F_{0}(s) \ \frac{J_{\alpha}(rs)}{(rs)^{\alpha}} \ s^{2 \alpha +1} \ ds,
\end{eqnarray}
where $J_{\alpha}$ denotes the Bessel function of the first kind.
\par Thus, in terms of the spherical mean operator (\ref{MF}) and according to (\ref{F01}) and (\ref{F02}) we can write (\ref{Uk}) as 
\begin{eqnarray}
\nonumber u_{k}(x,t) &=& \int_{\mathbb{R}^{n}} P_{k,t}^{11}(y) \ \tau_{x}(k) \ f(y) \ v_{k}(y) \ dy \\
\nonumber &+& \int_{\mathbb{R}^{n}} P_{k,t}^{12}(y) \ \tau_{x}(k) \ g(y) \ v_{k}(y) \ dy \\
\nonumber &=& \int_{0}^{+\infty} r^{2\gamma_{k}+n-1} \ \int_{S^{n-1}}P_{k,t}^{11}(ry') \ \tau_{x}(k) \ f(ry') \ v_{k}(y') \ dw (y') \ dr \\
\nonumber &+& \int_{0}^{+\infty} r^{2\gamma_{k}+n-1} \ \int_{S^{n-1}}P_{k,t}^{12}(ry') \ \tau_{x}(k) \ g(ry') \ v_{k}(y') \ dw (y') \ dr \\ 
\nonumber &=& d_{k} \ \int_{0}^{+\infty} r^{2\gamma_{k}+n-1} \ \mathcal{H}_{\gamma_{k}+n/2-1} \ F_{t}(r) \ \mathcal{M}_{f}(x,r) \ dr \\
 &+& d_{k} \ \int_{0}^{+\infty} r^{2\gamma_{k}+n-1} \ \mathcal{H}_{\gamma_{k}+n/2-1} \ G_{t}(r) \ \mathcal{M}_{g}(x,r) \ dr,
\end{eqnarray}
where 
\begin{eqnarray}\label{Fts}
 F_{t}(s) = \cos \big(t (\vert s \vert^{2} + m^{2})^{1/2}\big) 
\end{eqnarray} 
 and \\
 \begin{eqnarray}\label{Gts}
 G_{t}(s) = \frac{\sin \big(t (\vert s \vert^{2} + m^{2})^{1/2}\big)}{(\vert s \vert^{2} + m^{2})^{1/2}}  
\end{eqnarray}
\par On the other hand , we have on the basing of (\ref{F02}) and  (\ref{Gts})
\begin{eqnarray}\label{HGT}
 \mathcal{H}_{\alpha}G_{t}(r) = \frac{1}{2^{\alpha}\Gamma (\alpha +1)r^{\alpha}} \int_{0}^{+\infty} \frac{\sin (t(s^{2}+m^{2})^{1/2})}{(s^{2}+m^{2})^{1/2}} \ J_{\alpha} (rs) \ s^{\alpha +1}  ds 
 \end{eqnarray}

Later we shall determinated the value of integral
\begin{eqnarray}\label{intI}
I = \int_{0}^{+\infty} \frac{\sin (t(s^{2}+m^{2})^{1/2})}{(s^{2}+m^{2})^{1/2}} \ J_{\alpha} (rs) \ s^{\alpha +1} \ ds
\end{eqnarray}
Erdeyli has introduced may Hankel transform of many functions in \cite{Wat2} and   has defined the Hankel transform $\mathcal{H}_{\alpha}$ by
\begin{eqnarray}
\nonumber \mathcal{H}_{\alpha}(f)(y)= \int_{0}^{+\infty} f(x) J_{\nu}(xy) \ (xy)^{1/2} \ dx
\end{eqnarray} 
Using the formula $18$ in $($\cite{Wat2}, \ p. \ $59)$ \ and on the basing of (\ref{HGT}) we shall proceeding as follows: \\
 For \ \ $f(x) = x^{\nu +1/2} (x^{2}+\beta^{2})^{-1/2 \mu } \ J_{\mu} \big[a (x^{2}+\beta^{2})^{1/2}\big]$ \ \ such that \ $a > 0 , \ \Re (\beta) > 0$ \ and \ $Re (\mu) > \Re (\nu) > -1$ , 
we apply the Hankel transform  for  $y > 0$ to obtain
\begin{eqnarray}
\mathcal{H}_{\alpha}(f)(y) = a^{-\mu} y^{\nu +1/2} \beta^{-\mu+\nu+1} (a^{2}-y^{2})^{1/2\mu-1/2\nu-1/2} J_{\mu-\nu-1}[\beta (a^{2}-y^{2})^{1/2}]
\end{eqnarray}
Now , using the fact that \ $J_{1/2}(t)\ = \ (2/\pi)^{1/2} \sin (t) / \sqrt{t}$ \ and if we take in \ (\ref{Gts}) \ these values as:
$$\beta = m, \ \ \ \ \ \ \ \ \ \nu = 1/2, \ \ \ \ \ \ \ \ \  a = t, \ \ \ \ \ \ \ \ \ v = \alpha,   $$
 we obtain
\begin{eqnarray}\label{form1}
\nonumber  &\mathcal{H}_{\alpha}& \bigg( x^{\alpha+1/2} (x^{2}+m^{2})^{-1/4}   J_{1/2} \big[t (m^{2}+x^{2})^{1/2}\big]\bigg)(y) \\
\nonumber &=& \mathcal{H}_{\alpha} \bigg(\big(\frac{2}{\pi}\big)^{1/2} t^{-1/4} x^{\alpha +1/2} \frac{\sin (t (x^{2}+m^{2})^{1/2})}{(x^{2}+m^{2})^{1/2}}\bigg)(y) \\
\nonumber &=& \big(\frac{2 y}{\pi t}\big)^{1/2} \int_{0}^{+\infty}\frac{\sin (t (x^{2}+m^{2})^{1/2})}{(x^{2}+m^{2})^{1/2}} J_{\alpha}(x y) x^{\alpha+1} \ dx \\
&=& \big(\frac{2 y}{\pi t}\big)^{1/2} \ \ I
\end{eqnarray}
where $I$ is the integral defined in (\ref{intI}).\\
On the other hand , using the formula $(18)$ given in $($\cite{Wat2}, \ p. \ $59)$ \ we obtain
\begin{eqnarray}\label{form2}
\nonumber  &\mathcal{H}_{\alpha}& \bigg( x^{\alpha+1/2} (x^{2}+m^{2})^{-1/4}   J_{1/2} \big[t (m^{2}+x^{2})^{1/2}\big]\bigg)(y) \\
  &=& t^{-1/2} \big(\frac{y}{m}\big)^{\alpha+1/2} (t^{2}-y^{2})^{-1/2 \alpha - 1/2} J_{-\alpha-1/2}[m(t^{2}-y^{2})^{1/2}] \ \mathbb{1}_{\{0 < y < t\}} 
\end{eqnarray}
Now, according to (\ref{intI}), (\ref{form1})  and (\ref{form2}) we obtain
\begin{eqnarray}\label{result3}
\nonumber \mathcal{H}_{\alpha}G_{t}(r) &=& \frac{1}{2^{\alpha}\Gamma (\alpha +1)r^{\alpha}} \int_{0}^{+\infty} \frac{\sin (t(s^{2}+m^{2})^{1/2})}{(s^{2}+m^{2})^{1/2}} \ J_{\alpha} (rs) \ s^{\alpha +1} \ ds \\
\nonumber &=&  \left\{
    \begin{array}{ll}
      \frac{\sqrt{\pi}}{\Gamma (\alpha +1)} \big(\frac{m}{2}\big)^{\alpha +1/2} \ (t^{2}-y^{2})^{-\frac{1}{2} \alpha - 1/4} J_{-\alpha-1/2}[m(t^{2}-y^{2})^{1/2}]  , & \hbox{ if $0<y<t$},\\
0, & \hbox{ if $t<y<\infty$}
       \end{array}
  \right.
  \\
  &=& \frac{\sqrt{\pi}}{\Gamma (\alpha +1)} \big(\frac{m}{\sqrt{2}}\big)^{2 \alpha +1} \ \mathbb{S}_{-\alpha - 1/2}[m(t^{2}-y^{2})^{1/2}] 
\end{eqnarray}
Similarly, using the fact that \ $J_{-1/2}(t)\ = \ (2/\pi)^{1/2} \cos (t) $ \ and if we take in (\ref{Fts}) these values as :
$$\beta = m, \ \ \ \ \ \ \ \ \ \nu = -1/2, \ \ \ \ \ \ \ \ \  a = t, \ \ \ \ \ \ \ \ \ v = \alpha, $$
 we obtain
\begin{eqnarray}\label{result1}
\nonumber \mathcal{H}_{\alpha}F_{t}(r) &=& \frac{1}{2^{\alpha}\Gamma (\alpha +1)r^{\alpha}} \int_{0}^{+\infty} \cos (t(s^{2}+m^{2})^{1/2}) \ J_{\alpha} (rs) \ s^{\alpha +1} \ ds \\
\nonumber &=& \left\{
    \begin{array}{ll}
     \frac{\sqrt{\pi}}{\Gamma (\alpha +1)} 2^{-\alpha - 1/2} m^{\alpha +3/2}  \ t \ (t^{2}-y^{2})^{-\frac{1}{2} \alpha - 3/4} J_{-\alpha-3/2}[m(t^{2}-y^{2})^{1/2}] , & \hbox{ if $0<y<t$},\\
0, & \hbox{ if $t<y<\infty$}
    \end{array}
  \right.
  \\
  &=&  \frac{\sqrt{\pi}}{\Gamma (\alpha +1)} 2^{-\alpha - 1/2} m^{2 \alpha +3}  \ t \ \mathbb{S}_{-\alpha - 3/2}[m(t^{2}-y^{2})^{1/2}]
\end{eqnarray}
On the other hand, we have the recursion relation
\begin{eqnarray}\label{recursion}
\frac{d}{dx} \big[x^{\beta} \ J_{\beta} (x)\big] \ = \ x^{\beta} \ J_{\beta -1}(x)
\end{eqnarray}
It's directly by (\ref{recursion}) that
\begin{eqnarray}
\nonumber \frac{d}{dx} \mathbb{S}^{+}_{\beta} (x) \ = \ x \ \mathbb{S}^{+}_{\beta -1}(x)
\end{eqnarray}
where  $\mathbb{S}^{+}_{\beta}$  is the function defined in (\ref{Slmd}). Indeed
\begin{eqnarray}\label{result2}
\nonumber \frac{d}{dt} \big(\mathbb{S}_{-\alpha-1/2} [m (t^{2}-r^{2})^{1/2}]\big)  &=& \frac{d}{dt} [m (t^{2}-r^{2})^{1/2}] \ \big(\frac{d}{dt} \mathbb{S}_{-\alpha-1/2}\big)[m (t^{2}-r^{2})^{1/2}] \\
\nonumber &=& m^{2} \ t \ \mathbb{S}_{-\alpha-3/2}[m (t^{2}-r^{2})^{1/2}]
\end{eqnarray}
Then, we obtain
\begin{eqnarray}
\mathbb{S}_{-\alpha-3/2}^{+}[m (t^{2}-r^{2})^{1/2}] = \frac{1}{m^{2} t} \ \frac{d}{dt} \big(\mathbb{S}_{-\alpha-1/2}^{+} [m (t^{2}-r^{2})^{1/2}]\big)
\end{eqnarray}
According to (\ref{result1}) and (\ref{result2}) we deduce that
\begin{eqnarray}\label{result4}
\nonumber \mathcal{H}_{\alpha}F_{t}(r) &=& \frac{1}{2^{\alpha}\Gamma (\alpha +1)r^{\alpha}} \int_{0}^{+\infty} \cos (t(s^{2}+m^{2})^{1/2}) \ J_{\alpha} (rs) \ s^{\alpha +1} \ ds \\
  &=& \frac{\sqrt{\pi}}{\Gamma (\alpha +1)} \big(\frac{m}{\sqrt{2}}\big)^{2 \alpha +1}  \ \frac{d}{dt} \big(\mathbb{S}_{-\alpha-1/2}^{+} [m (t^{2}-r^{2})^{1/2}]\big)
\end{eqnarray}
 Consequently, according to (\ref{Uk}), (\ref{result3}) and (\ref{result4}), the theorem \ref{thp} holds.
\par Now, we are in position to prove the lemma \ref{ESTM}.

\begin{proof}(of lemma \ref{ESTM}) \\
In \ view \ of \ the \ Dunkl's solution \ of \ the \ Klein-Gordon-equation \  given \ by \ (\ref{Sol}) we obtain by \ mean \ of \ the \ Plancherel \ formula
\begin{eqnarray}
\nonumber && \int_{\mathbb{R}^{n}}   \vert u_{k}(x,t)\vert^{2}  \ w_{k}(\xi) \ d\xi  \ = \ \int_{\mathbb{R}^{n}}\vert \mathcal{F}_{k} (u_{k}(. , t)(\xi) \vert^{2} \ w_{k}(\xi) \ d\xi \\ 
\nonumber &=& \int_{\mathbb{R}^{n}} \mathcal{F}_{k} (u_{k}(. , t)(\xi) \ \overline{\mathcal{F}_{k} (u_{k}(. , t)(\xi)} \ w_{k}(x) \ dx \\
\nonumber &=& \int_{\mathbb{R}^{n}} \bigg[\frac{\sin(t (\vert \xi \vert^{2}+m^{2})^{1/2})}{(\vert \xi \vert^{2}+m^{2})^{1/2}} \ \mathcal{F}_{k}(f)(\xi)  \ + \ \cos (t (\vert \xi \vert^{2}+m^{2})^{1/2}) \mathcal{F}_{k}(g)(\xi) \bigg] \\ 
\nonumber &\times & \bigg[\frac{\sin(t (\vert \xi \vert^{2}+m^{2})^{1/2})}{(\vert \xi \vert^{2}+m^{2})^{1/2}} \ \overline{\mathcal{F}_{k}(f)(\xi)} \ + \ \cos (t (\vert \xi \vert^{2}+m^{2})^{1/2}) \ \overline{\mathcal{F}_{k}(g)(\xi)}\bigg] \\
\nonumber &=&  \int_{\mathbb{R}^{n}}\frac{\sin^{2}(t (\vert \xi \vert^{2}+m^{2})^{1/2})}{(\vert \xi \vert^{2}+m^{2})^{1/2}} \ \mathcal{F}_{k}(f)(\xi) \overline{\mathcal{F}_{k}(f)(\xi)} \  \ w_{k}(\xi) \ d\xi \\ 
\nonumber &+&  \int_{\mathbb{R}^{n}}\frac{\sin(t (\vert \xi \vert^{2}+m^{2})^{1/2})}{(\vert \xi \vert^{2}+m^{2})^{1/2}} \ \cos (t (\vert \xi \vert^{2}+m^{2})^{1/2}) \ \mathcal{F}_{k}(f)(\xi) \overline{\mathcal{F}_{k}(g)(\xi)} \  \ w_{k}(\xi) \ d\xi \\
\nonumber &+&  \int_{\mathbb{R}^{n}}\frac{\cos (t (\vert \xi \vert^{2}+m^{2})^{1/2})}{(\vert \xi \vert^{2}+m^{2})^{1/2}} \  \sin (t (\vert \xi \vert^{2}+m^{2})^{1/2}) \ \mathcal{F}_{k}(g)(\xi) \overline{\mathcal{F}_{k}(f)(\xi)} \  \ w_{k}(\xi) \ d\xi \\
\nonumber &+&  \int_{\mathbb{R}^{n}} \cos^{2} (t (\vert \xi \vert^{2}+m^{2})) \ \mathcal{F}_{k}(g)(\xi) \overline{\mathcal{F}_{k}(g)(\xi)} \  \ w_{k}(\xi) \ d\xi \\
\nonumber &=&  \int_{\mathbb{R}^{n}}\frac{\sin^{2} (t (\vert \xi \vert^{2}+m^{2})^{1/2})}{(\vert \xi \vert^{2}+m^{2})} \  \big\vert \mathcal{F}_{k}(f)(\xi)\big\vert^{2} \  \ w_{k}(\xi) \ d\xi \\
\nonumber &+& \cos^{2} (t (\vert \xi \vert^{2}+m^{2})^{1/2}) \  \big\vert \mathcal{F}_{k}(g)(\xi)\big\vert^{2} \  \ w_{k}(\xi) \ d\xi \\
\nonumber &+& \frac{1}{2} \int_{\mathbb{R}^{n}}  \frac{\sin (2 t (\vert \xi \vert^{2}+m^{2})^{1/2})}{(\vert \xi \vert^{2}+m^{2})^{1/2}} \bigg(\mathcal{F}_{k}(f)(\xi) \overline{\mathcal{F}_{k}(g)(\xi)} \ + \  \overline{\mathcal{F}_{k}(f)(\xi)} \mathcal{F}_{k}(g)(\xi)\bigg) \ w_{k}(\xi) \ d\xi 
\end{eqnarray}
Now, using these formulas:
\begin{eqnarray}\label{TRIG}
 \cos^{2}(x) \ = \ \frac{1 \ + \ \cos (2x)}{2}  \ \ \ \ and \ \ \ \ \sin^{2}(x) \ = \ \frac{1 \ - \ \cos (2x)}{2} \ \ \ \ \forall \ \ x \in \mathbb{R}
\end{eqnarray}
we obtain
\begin{eqnarray}\label{I1I2I3}
\nonumber && \int_{\mathbb{R}^{n}}   \vert u_{k}(x,t)\vert^{2}  \ w_{k}(\xi) \ d\xi  \ = \ \frac{c_{k}^{-2}}{2} \int_{\mathbb{R}^{n}} \bigg\{\frac{\vert \mathcal{F}_{k}(f)(\xi)\vert^{2}}{(\vert \xi \vert^{2}+m^{2})} \ + \ \vert \mathcal{F}_{k}(g)(\xi)\vert^{2} \bigg\}  \ w_{k}(\xi) \ d\xi \\
\nonumber &-&  \frac{c_{k}^{-2}}{2} \int_{\mathbb{R}^{n}} \frac{\vert \mathcal{F}_{k}(f)(\xi)\vert^{2}}{(\vert \xi \vert^{2}+m^{2})} \cos (2 t (\vert \xi \vert^{2}+m^{2})^{1/2}) \ w_{k}(\xi) \ d\xi \\
\nonumber &+&  \frac{c_{k}^{-2}}{2} \int_{\mathbb{R}^{n}}  \cos (2 t (\vert \xi \vert^{2}+m^{2})^{1/2}) \  \vert \mathcal{F}_{k}(g)(\xi)\vert^{2}  \ w_{k}(\xi)\ d\xi \\
\nonumber &+&  \frac{c_{k}^{-2}}{2} \int_{\mathbb{R}^{n}}\frac{\sin (2 t (\vert \xi \vert^{2}+m^{2})^{1/2})}{(\vert \xi \vert^{2}+m^{2})^{1/2}} \bigg(\mathcal{F}_{k}(f)(\xi) \overline{\mathcal{F}_{k}(g)(\xi)} \ + \ \overline{\mathcal{F}_{k}(f)(\xi)} \mathcal{F}_{k}(g)(\xi) \bigg)  \ w_{k}(\xi)\ d\xi \\
\nonumber &=& \frac{c_{k}^{-2}}{2} \int_{\mathbb{R}^{n}} \bigg\{\frac{\vert \mathcal{F}_{k}(f)(\xi)\vert^{2}}{(\vert \xi \vert^{2}+m^{2})} \ + \ \vert \mathcal{F}_{k}(g)(\xi)\vert^{2} \bigg\}  \ w_{k}(\xi) \ d\xi \\
 &-& I_{1} + I_{2} + I_{3}
\end{eqnarray}
where
\begin{eqnarray}
\nonumber I_{1} &=& \frac{c_{k}^{-2}}{2} \int_{\mathbb{R}^{n}} \frac{\vert \mathcal{F}_{k}(f)(\xi)\vert^{2}}{(\vert \xi \vert^{2}+m^{2})} \cos (2 t (\vert \xi \vert^{2}+m^{2})^{1/2}) \ w_{k}(\xi) \ d\xi \\
\nonumber I_{2} &=& \frac{c_{k}^{-2}}{2} \int_{\mathbb{R}^{n}}  \cos (2 t (\vert \xi \vert^{2}+m^{2})^{1/2}) \  \vert \mathcal{F}_{k}(g)(\xi)\vert^{2}  \ w_{k}(\xi)\ d\xi \\
\nonumber I_{3} &=& \frac{c_{k}^{-2}}{2} \int_{\mathbb{R}^{n}}\frac{\sin (2 t (\vert \xi \vert^{2}+m^{2})^{1/2})}{(\vert \xi \vert^{2}+m^{2})^{1/2}} \bigg(\mathcal{F}_{k}(f)(\xi) \overline{\mathcal{F}_{k}(g)(\xi)} \ + \ \overline{\mathcal{F}_{k}(f)(\xi)} \mathcal{F}_{k}(g)(\xi) \bigg)  \ w_{k}(\xi)\ d\xi 
\end{eqnarray}
It is easy by the Riemann-Lebesgue theorem to prove that
\begin{eqnarray}\label{limit}
\lim_{\vert t \vert \longmapsto + \infty} I_{1} \ = \ \lim_{\vert t \vert \longmapsto + \infty} I_{2} \ = \ \lim_{\vert t \vert \longmapsto + \infty} I_{3} \ = \ 0
\end{eqnarray}
According to (\ref{I1I2I3}) and (\ref{limit}), we deduce that
\begin{eqnarray}
\nonumber \lim_{\vert t \vert \longmapsto + \infty} \int_{\mathbb{R}^{n}}
\vert u_{k} (x , t) \vert^{2} \ w_{k}(x) \ dx &=& \frac{c_{k}^{-2}}{2} \int_{\mathbb{R}^{n}} \bigg\{ \vert \mathcal{F}_{k} g (\xi) \vert^{2} + \frac{\vert \mathcal{F}_{k} f (\xi) \vert^{2}}{(\vert \xi \vert^{2}+m^{2})} \bigg\} \ w_{k}(x) \ dx  \\
\nonumber &=& \frac{1}{2} \Vert g \Vert_{k}^{2} \ + \ \frac{1}{2} \bigg\Vert \frac{\mathcal{F}_{k}(f)(.)}{\vert . \vert^{2} + m^ {2}} \bigg\Vert_{k}^{2} \\
\nonumber &=& \frac{1}{2} \Vert g \Vert_{k}^{2} \ + \ \frac{1}{2} \big\Vert \big(- \Delta_{k} + m^{2} \big)^{-1/2} f \big\Vert_{k}^{2}.
\end{eqnarray}
We have 
\begin{eqnarray}
\nonumber \int_{\mathbb{R}^{n}}\frac{\vert \mathcal{F}_{k}(f)(\xi)(\xi) \vert^{2}}{(\vert \xi \vert^{2} \ + \ m^{2})} \ w_{k}(\xi)\ d\xi \ \leqslant \ \int_{\mathbb{R}^{n}}\frac{\vert \mathcal{F}_{k}(f)(\xi)(\xi) \vert^{2}}{\vert \xi \vert^{2}} \ w_{k}(\xi)\ d\xi \ = \ \big\Vert \big(-\Delta_{k}\big)^{-1/2}f \big\Vert_{k}^{2}.
\end{eqnarray}
Then
 \begin{eqnarray}\label{Strichartz}
  \lim_{\vert t \vert \longmapsto + \infty } \Vert u_{k} (. , t) \Vert_{k}^{2}
  \leqslant \ \frac{1}{2} \ \Vert (- \Delta_{k})^{-1/2} f \Vert_{k}^{2} + \frac{1}{2} \ \Vert g\Vert_{k}^{2}
 \end{eqnarray}
 This completes the proof of lemma \ref{ESTM} and the estimate (\ref{estm}) holds. 
 \end{proof}
\textbf{Comments:}\\
The estimate (\ref{Strichartz}) have big importance in physical . We restrict our attention to the $L^{2}$- behaviors because these are the most physically interesting quantities . \\
Indeed, we have proving the following Strichartz-type inequality:
\begin{eqnarray}
\nonumber \Vert u_{k} (. , t) \Vert_{k}^{2} \ \leqslant \  \Vert (- \Delta_{k})^{-1/2} f \Vert_{k}^{2} + \ \Vert g\Vert_{k}^{2} 
 \end{eqnarray}
with \ $\Vert.\Vert_{k}$ \ denotes the norm in \ $L^{2} (\mathbb{R}^{n} ; w_{k}(x) dx)$. \\
Consequently, it follows that if \ \ $\Vert u_{k} (. , t) \Vert_{k}^{2} \rightarrow 0$ \ as \ $\vert t \vert \rightarrow +\infty$ \ then \ \ $u_{k} \equiv 0$.

\section{Energy \ Associated \ To \ The \ Dunkl \ Klein-Gordon \ Equation}\label{sec5}
In this section we only assume $k \in \mathcal{R}^{+}$ and $n \geqslant 1$. Let $u_{k}(x,t)$ the solution of the Klein-Gordon equation. Define the kinetic and potential energies by:
\begin{eqnarray}\label{kin}
\mathcal{K}_{k}[u_{k}](t) &=& \frac{1}{2}\int_{\mathbb{R}^{n}}\vert \partial_{t} u_{k}(x,t)^{2} \ w_{k}(x) \ dx 
\end{eqnarray}
and
\begin{eqnarray}\label{pot}
\mathcal{P}_{k}[u_{k}](t) &=& \frac{1}{2} \int_{\mathbb{R}^{n}}\sum_{j=1}^{n}\vert T_{j}^{x} (k)(u_{k}(x;t))\vert^{2} \ w_{k}(x) \ dx 
\end{eqnarray}
where \ $T_{j}^{x} (k)$ \ is the Dunkl operator. \\
By definition, the total energy of $u_{k}$ is 
\begin{eqnarray}\label{total}
\mathcal{E}_{k}[u_{k}](t) = \mathcal{K}_{k}[u_{k}](t) + \mathcal{P}_{k}[u_{k}](t)
\end{eqnarray}
Since
\begin{eqnarray}
\nonumber \mathcal{F}_{k} \big(T_{j}^{x}(k) u_{k}(. , t)\big)(\xi) = - i \xi_{j} \ \mathcal{F}_{k} (u_{k}(. , t)(\xi),
\end{eqnarray}
then according to (\ref{pot}) and (\ref{total}), we obtain by means of the Plancherel formula  
\begin{eqnarray}
 \nonumber \mathcal{E}_{k}[u_{k}](t) = \frac{c_{k}^{-2}}{2} \int_{\mathbb{R}^{n}} \bigg\{ \big\vert \partial_{t} \mathcal{F}_{k} (u_{k}(. , t)(\xi) \big\vert^{2} + \vert \xi \vert^{2} \ \big\vert  \mathcal{F}_{k} (u_{k}(. , t)(\xi) \big\vert^{2} \bigg\} \ w_{k}(\xi) \ d \xi
\end{eqnarray}
On the other hand , using (\ref{Sol}) we can write
\begin{eqnarray}
\nonumber \big\vert  \mathcal{F}_{k} (u_{k}(. , t)(\xi) \big\vert^{2} &=& \sin^{2} (t(\vert \xi \vert^{2} + m^{2})^{1/2}) \ \big\vert  \mathcal{F}_{k} (f)(\xi) \big\vert^{2} \\
\nonumber \nonumber &+&  \cos^{2} (t(\vert \xi \vert^{2} + m^{2})^{1/2}) \ \big\vert  \mathcal{F}_{k} (g)(\xi) \big\vert^{2} \\
\nonumber &+&  \ \frac{\sin (2 t(\vert \xi \vert^{2} + m^{2})^{1/2})}{(\vert \xi \vert^{2} + m^{2})^{1/2}} \ \Re \big(\mathcal{F}_{k} (f)(\xi) \ \overline{\mathcal{F}_{k} (g)(\xi)}\big)
\end{eqnarray}
Similarly, we obtain
\begin{eqnarray}\label{partiall}
\nonumber  &\big\vert   \partial_{t}  \mathcal{F}_{k} (u_{k}(. , t)(\xi) \big\vert^{2}& = \cos^{2} (t(\vert \xi \vert^{2} + m^{2})^{1/2}) \ \big\vert  \mathcal{F}_{k} (f)(\xi) \big\vert^{2} \\
\nonumber &+& (\vert \xi \vert^{2} + m^{2}) \ \sin^{2} (t(\vert \xi \vert^{2} + m^{2})^{1/2}) \ \big\vert  \mathcal{F}_{k} (g)(\xi) \big\vert^{2} \\
\nonumber &-&  \ (\vert \xi \vert^{2} + m^{2})^{1/2} \ \sin (2 t(\vert \xi \vert^{2} + m^{2})^{1/2}) \  \Re \big(\mathcal{F}_{k} (f)(\xi) \ \overline{\mathcal{F}_{k} (g)(\xi)}\big) \\
\end{eqnarray}
Thus on the based of (\ref{kin}) and (\ref{partiall}) we obtain
\begin{eqnarray}\label{kiok}
\nonumber & \mathcal{K}_{k}[u_{k}](t)& = \frac{c_{k}^{-2}}{2} \int_{\mathbb{R}^{n}}\cos^{2} (t(\vert \xi \vert^{2} + m^{2})^{1/2}) \ \big\vert  \mathcal{F}_{k} (f)(\xi) \big\vert^{2} \ w_{k}(\xi) d\xi \\
\nonumber &+&  \frac{c_{k}^{-2}}{2} \int_{\mathbb{R}^{n}}(\vert \xi \vert^{2} + m^{2}) \ \sin^{2} (t(\vert \xi \vert^{2} + m^{2})^{1/2}) \ \big\vert  \mathcal{F}_{k} (g)(\xi) \big\vert^{2} \ w_{k}(\xi) d\xi \\
\nonumber &-& \frac{c_{k}^{-2}}{2} \int_{\mathbb{R}^{n}}  (\vert \xi \vert^{2} + m^{2})^{1/2})\sin (2 t(\vert \xi \vert^{2} + m^{2})^{1/2})  \ \Re \big(\mathcal{F}_{k} (f)(\xi) \ \overline{\mathcal{F}_{k} (g)(\xi)}\big)\ w_{k}(\xi) d\xi \\
\end{eqnarray}
Now, using (\ref{TRIG}) we can write (\ref{kiok}) as
\begin{eqnarray}\label{bnv}
\nonumber & \mathcal{K}_{k}[u_{k}](t)& = \frac{c_{k}^{-2}}{4} \int_{\mathbb{R}^{n}} \big\vert  \mathcal{F}_{k} (f)(\xi) \big\vert^{2} \ w_{k}(\xi) d\xi \\
\nonumber &+&  \frac{c_{k}^{-2}}{4} \int_{\mathbb{R}^{n}} \cos (2 t (\vert \xi \vert^{2} + m^{2})^{1/2}) \ \big\vert  \mathcal{F}_{k} (f)(\xi) \big\vert^{2} \  w_{k}(\xi) d\xi \\
\nonumber &+&   \frac{c_{k}^{-2}}{4} \int_{\mathbb{R}^{n}} (\vert \xi \vert^{2} + m^{2}) \ \big\vert  \mathcal{F}_{k} (g)(\xi) \big\vert^{2} \ w_{k}(\xi) d\xi \\
\nonumber &-&   \frac{c_{k}^{-2}}{4} \int_{\mathbb{R}^{n}} (\vert \xi \vert^{2} + m^{2}) \ \cos (2 t (\vert \xi \vert^{2} + m^{2})^{1/2}) \ \big\vert  \mathcal{F}_{k} (g)(\xi) \big\vert^{2} \ w_{k}(\xi) d\xi \\
\nonumber &-&   \frac{c_{k}^{-2}}{4} \ \int_{\mathbb{R}^{n}} (\vert \xi \vert^{2} + m^{2})^{1/2}) \ \sin (2 t (\vert \xi \vert^{2} + m^{2})^{1/2}) \ \Re \big(\mathcal{F}_{k} (f)(\xi) \overline{\mathcal{F}_{k} (g)(\xi)}\big) \ w_{k}(\xi) d\xi \\
\nonumber &=& \frac{c_{k}^{-2}}{4} \Vert \mathcal{F}_{k} (f) \Vert_{k}^{2} + \frac{c_{k}^{-2}}{4} \Vert \langle. , .\rangle^{1/2}  \mathcal{F}_{k} (g) \Vert_{k}^{2} +  \frac{c_{k}^{-2}}{4} m^{2} \Vert \mathcal{F}_{k} (g) \Vert_{k}^{2} \\
\nonumber &+&  \frac{c_{k}^{-2}}{4} \int_{\mathbb{R}^{n}} \big[  \big\vert  \mathcal{F}_{k} (f)(\xi) \big\vert^{2} - (\vert \xi \vert^{2} + m^{2})  \big\vert  \mathcal{F}_{k} (g)(\xi) \big\vert^{2} \big]  \cos (2 t (\vert \xi \vert^{2} + m^{2})^{1/2}) \  w_{k}(\xi) d\xi \\ 
 \nonumber &-&  \frac{c_{k}^{-2}}{4} \int_{\mathbb{R}^{n}} \big[\mathcal{F}_{k} (f)(\xi) \overline{\mathcal{F}_{k} (g)(\xi)} + \overline{\mathcal{F}_{k} (f)(\xi)} \mathcal{F}_{k} (g)(\xi) \big] \ \sin (2 t (\vert \xi \vert^{2} + m^{2})^{1/2}) \  w_{k}(\xi) d\xi \\
\end{eqnarray}

We apply the limit as \ $\vert t \vert \rightarrow +\infty$ \ to (\ref{bnv}) we obtain on the based of the Riemann-Lebesgue theorem
\begin{eqnarray}\label{limkin}
\nonumber \lim_{\vert t \vert \rightarrow +\infty} \mathcal{K}_{k}[u_{k}](t) \ &=& \  \frac{c_{k}^{-2}}{4} \int_{\mathbb{R}^{n}} \big\vert  \mathcal{F}_{k} (f)(\xi) \big\vert^{2} \ w_{k}(\xi) d\xi \\
 \nonumber &+&  \frac{c_{k}^{-2}}{4} \int_{\mathbb{R}^{n}} (\vert \xi \vert^{2} + m^{2}) \ \big\vert  \mathcal{F}_{k} (g)(\xi) \big\vert^{2} \ w_{k}(\xi) d\xi \\
 \nonumber &=&  \frac{c_{k}^{-2}}{4} \Vert \mathcal{F}_{k}(f) \Vert_{k}^{2} \ + \  \frac{c_{k}^{-2}}{4} \Vert \langle . , . \rangle^{1/2} \mathcal{F}_{k}(g) \Vert_{k}^{2} \ + \  \frac{c_{k}^{-2}}{4} m^{2}  \Vert \mathcal{F}_{k}(g) \Vert_{k}^{2} \\
&=&  \frac{c_{k}^{-2}}{4} \Vert f \Vert_{k}^{2} \ + \  \frac{c_{k}^{-2}}{4} \Vert T_{j}^{k}(g) \Vert_{k}^{2} \ + \  \frac{c_{k}^{-2}}{4} m^{2} \Vert g \Vert_{k}^{2}
\end{eqnarray}
On the other hand, we obtain according to (\ref{TRIG}) and (\ref{pot})
\begin{eqnarray}\label{bnp}
\nonumber &\mathcal{P}_{k}[u_{k}](t)& = \frac{c_{k}^{-2}}{4}  \int_{\mathbb{R}^{n}} \vert \xi \vert^{2} \ \vert \mathcal{F}_{k}(u_{k}(. , t)) \vert^{2} \ w_{k}(\xi) d\xi \\
\nonumber &=& \frac{c_{k}^{-2}}{4} \int_{\mathbb{R}^{n}} \frac{\vert \xi \vert^{2}}{\vert \xi \vert^{2} + m^{2}} \vert \mathcal{F}_{k}(f)(\xi)  \vert^{2} \ w_{k}(\xi) d\xi \\
\nonumber &-& \frac{c_{k}^{-2}}{4}  \int_{\mathbb{R}^{n}} \frac{\vert \xi \vert^{2}}{\vert \xi \vert^{2} + m^{2}} \cos (2 t(\vert \xi \vert^{2} + m^{2})^{1/2}) \vert \mathcal{F}_{k}(f)(\xi)  \vert^{2} \ w_{k}(\xi) d\xi \\
\nonumber &+& \frac{c_{k}^{-2}}{4}  \int_{\mathbb{R}^{n}} \vert \xi \vert^{2} \ \vert \mathcal{F}_{k}(g)(\xi)  \vert^{2} \ w_{k}(\xi) d\xi \\
\nonumber &+& \frac{c_{k}^{-2}}{4}  \int_{\mathbb{R}^{n}} \vert \xi \vert^{2} \ \cos (2 t(\vert \xi \vert^{2} + m^{2})^{1/2}) \vert \mathcal{F}_{k}(g)(\xi)  \vert^{2} \ w_{k}(\xi) d\xi \\
\nonumber &+& \frac{c_{k}^{-2}}{4}  \int_{\mathbb{R}^{n}}  \frac{\vert \xi \vert^{2} \sin (2 t(\vert \xi \vert^{2} + m^{2})^{1/2})}{(\vert \xi \vert^{2} + m^{2})}  \big(\mathcal{F}_{k}(f)(\xi)  \overline{\mathcal{F}_{k}(g)(\xi)} +  \ \overline{\mathcal{F}_{k}(f)(\xi)}\mathcal{F}_{k}(g)(\xi)\big)  w_{k}(\xi) d\xi \\
\nonumber &=&  \frac{c_{k}^{-2}}{4} \Vert \langle . , . \rangle^{1/2} \mathcal{F}_{k}(g) \Vert_{k}^{2} + \frac{c_{k}^{-2}}{4} \Vert \mathcal{F}_{k}(f) \Vert_{k}^{2} - \frac{c_{k}^{-2}}{4} m^{2}  \int_{\mathbb{R}^{n}} \frac{1}{(\vert \xi \vert^{2} + m^{2})} \Vert \mathcal{F}_{k}(f) \Vert_{k}^{2} \ w_{k}(\xi) d\xi  \\
\nonumber &+& \frac{c_{k}^{-2}}{4}  \int_{\mathbb{R}^{n}} \big[ \vert \xi \vert^{2} \vert \mathcal{F}_{k}(g)(\xi)  \vert^{2} -  \vert \mathcal{F}_{k}(f)(\xi)  \vert^{2} \big] \cos (2 t(\vert \xi \vert^{2} + m^{2})^{1/2})  \ w_{k}(\xi) d\xi  \\
\nonumber &+& \frac{c_{k}^{-2}}{4} m^{2} \int_{\mathbb{R}^{n}} \frac{1}{(\vert \xi \vert^{2} + m^{2})}\vert \mathcal{F}_{k}(f)(\xi)  \vert^{2}  \cos (2 t(\vert \xi \vert^{2} + m^{2})^{1/2})  \ w_{k}(\xi) d\xi \\
 \nonumber &+& \frac{c_{k}^{-2}}{4}  \int_{\mathbb{R}^{n}} \frac{\vert \xi \vert^{2} \ \sin (2 t (\vert \xi \vert^{2} + m^{2})^{1/2})}{(\vert \xi \vert^{2} + m^{2})^{1/2}}   \big[\mathcal{F}_{k} (f)(\xi) \overline{\mathcal{F}_{k} (g)(\xi)} + \overline{\mathcal{F}_{k} (f)(\xi)} \mathcal{F}_{k} (g)(\xi) \big]  \  w_{k}(\xi) d\xi \\
\end{eqnarray}

Using the  Riemann-Lebesgue theorem, we apply the limit to (\ref{bnp}) as  $\vert t \vert \rightarrow +\infty$  we obtain
\begin{eqnarray}\label{limpot}
 \nonumber \lim_{\vert t \vert \rightarrow +\infty} \mathcal{P}_{k}[u_{k}](t) \ &=&  \frac{c_{k}^{-2}}{4} \Vert f \Vert_{k}^{2} \ + \ \frac{c_{k}^{-2}}{4} \Vert T_{j}^{k}(g) \Vert_{k}^{2} \\
 \nonumber &-& \frac{c_{k}^{-2}}{4} m^{2}  \int_{\mathbb{R}^{n}} \frac{1}{(\vert \xi \vert^{2} + m^{2})}\Vert \mathcal{F}_{k}(f) \Vert_{k}^{2}  \  w_{k}(\xi) d\xi \\
 \nonumber &=& \frac{c_{k}^{-2}}{4} \Vert f \Vert_{k}^{2} \ + \ \frac{c_{k}^{-2}}{4} \Vert T_{j}^{k}(g) \Vert_{k}^{2} - \frac{c_{k}^{-2}}{4} m^{2} \big\Vert  \frac{1}{(\vert . \vert^{2} + m^{2})^{1/2}} \mathcal{F}_{k}(f)  \big\Vert_{k}^{2} \\
  &=& \frac{c_{k}^{-2}}{4} \Vert f \Vert_{k}^{2} \ + \ \frac{c_{k}^{-2}}{4} \Vert T_{j}^{k}(g) \Vert_{k}^{2} - \frac{c_{k}^{-2}}{4} m^{2} \Vert (-\Delta_{k}+m^{2})^{-1/2} f \Vert_{k}^{2} 
\end{eqnarray}
Consequently, on the based of (\ref{total}), (\ref{limkin}) and (\ref{limpot}) we deduce that
\begin{eqnarray}
 \nonumber \lim_{\vert t \vert \rightarrow +\infty} \mathcal{E}_{k}[u_{k}](t) &=& \frac{c_{k}^{-2}}{2} \Vert f \Vert_{k}^{2} \ + \ \frac{c_{k}^{-2}}{2} \Vert T_{j}^{k}(g) \Vert_{k}^{2} \ + \ \frac{m^{2} \ c_{k}^{-2}}{4} \Vert g \Vert_{k}^{2} \\
 \nonumber &-& \frac{c_{k}^{-2}}{4} m^{2} \Vert  \frac{1}{(\vert . \vert^{2} + m^{2})^{1/2}} \mathcal{F}_{k}(f)  \Vert_{k}^{2} \\
\nonumber &=& \frac{c_{k}^{-2}}{2} \Vert f \Vert_{k}^{2} \ + \ \frac{c_{k}^{-2}}{2} \Vert T_{j}^{k}(g) \Vert_{k}^{2} \ + \ \frac{m^{2} \ c_{k}^{-2}}{4} \Vert g \Vert_{k}^{2} \\
\nonumber &-& \frac{c_{k}^{-2}}{4} m^{2} \Vert (-\Delta_{k}+m^{2})^{-1/2} f \Vert_{k}^{2} 
\end{eqnarray}
Note that, the proposition \ref{collect} collects all the above facts and discussions. \\
\textbf{Comments:}\\


\begin{thebibliography}{00}


 \bibitem{AG2} B.  Amri and M. Gaidi,  \textit{$L^p-L^q$ estimates for the solution of the Dunkl wave equation}, manuscripta math. 159, 379-396 (2019).
\bibitem{AG} B.  Amri and M. Gaidi\textit{ $L^p$- Estimates for an Oscillating
Dunkl Multiplier}, Mediterr. J. Math. (2018) 15:85.
\bibitem{J1}M.F.E. de Jeu,\textit{The Dunkl transform},Invent. Math. 113 (1993), no. 1, 147-162.
\bibitem{D1} C. F. Dunkl, \textit{Differential-Difference operators associated to reflextion groups},
Trans. Amer. Math. 311 (1989), no. 1, 167--183.
\bibitem{D2} C. F. Dunkl, \textit{Hankel transforms associated to finite reflection groups}, Contemp. Math., vol. 138, 1992, pp. 123-138.
\bibitem{D3} C.F. Dunkl  \textit{Integral kernels with reflection group invariance}, Canad. J. Math. 43 (1991), 1213-1227.
\bibitem{E} Elcin Yusufoglu, \textit{The variational iteration method for studying the
Klein-Gordon equation}, Applied Mathematics Letters 21 (2008) 669-674.
\bibitem{Wat2} Erd\'eyli , A., \textit{et al.-Tables of Integral Transforms , Vol.II,}
McGrawHill,New York , 1954.
\bibitem{Feff} C. Feffermann and E. M. Stein,\textit{$H^p$ spaces of several variables}, Acta. Math. 129
(1972), 137-193.
\bibitem{H}   I. Hirschmann, \textit{On multiplier transformations},  Duke Math. J. 26 (1959), 221--242.
\bibitem{Horm}
L. H\"{o}rmander, \textit{Estimates for translation invariant operators in $L^p$ spaces}, Acta Math. 104
(1960), 93-140.
\bibitem{AD} W.Littman, \textit{The wave operator and $Lp$ norms}, Trans.Amer.Math.12(1963) , 55-68.
\bibitem{Majj1} H. Mejjaoli, \textit{Strichartz estimates for the Dunkl wave equation and application},
J.  Math. Anal. Appl. 346 (1) (2008), 41-54.
\bibitem{Majj2} H. Mejjaoli, \textit{Nonlinear generalized Dunkl-wave equations and applications}
 J. Math. Anal. Appl. 375 (2011) 118-138.
\bibitem{Voit} M. R$\ddot{o}$sler and M. Voit,  \textit{Markov processes with Dunkl operators}, Adv.in Appl.Math.21, (1988), 575, 643.
\bibitem{MIY} A. Miyachi,  \textit{ On some Fourier Multipliers for $H^p(\mathbb{R}^n))$}, J. Fac. Sci. Unvi. Tokyo, Sect IA, 27, 157-179.
\bibitem{Per}  C. Peral, \textit{$L^p $-estimates for the wave equation}, J. Funct. Anal., 36(1):114-145, 1980.
\bibitem{said1} S. B. Said and  B. \O rsted, \textit{The wave equation for Dunkl operators}, Indag. Math. (N.S.) 16 (2005), 351-391.
\bibitem{Stein2}
E. M. Stein, \textit{Interpolation of linear operators}, Transactions of
the American Mathematical Society 83 (1956), 482-492.
\bibitem{Sti} R. Strichartz,\textit{ Convolutions with kernels having singularities on a sphere}, Trans.Amer. Math. Soc. 148 (1970), 461-471.
\bibitem{Tr1} K. Trim\`eche,\textit{The Dunkl intertwining operator on spaces of functions and distributions and integral representation of its dual}, Integral Transform. Spec. Funct. 12 (2001), 349-374.
\bibitem{Tr2} K. Trim\`eche,\textit{Paley-Wiener theorems for the Dunkl transform and Dunkl translation operators}, Integral Transform. Spec. Funct. 13 (2002), 17-38.
\bibitem{Wat1} G.N. Watson, \textit{A Treatise on the Theory of Bessel Functions}
(2nd ed.). Cambridge Univ. Press, 1944.

\end{thebibliography}
\end{document}